\let\OToneAccents=\relax
\magnification=\magstep1
\frenchspacing
\baselineskip=16truept
\font\douze=cmr10 at 12pt

\font\grande=cmb10 at 16truept
\def\titre#1{{\OToneAccents\noindent\grande #1}}

\def \chapter#1{\vfill\eject\ifodd\pageno\else\ \vfill\eject\fi\centerline{\grande #1}\bigskip}

\font\tendo=wncyr10 at 12pt
\font\sevendo=wncyr7
\newfam\dofam 
\textfont\dofam = \tendo 
\scriptfont\dofam= \sevendo

\font\tendb=msbm10 
\font\sevendb=msbm7
\newfam\dbfam 
\textfont\dbfam = \tendb 
\scriptfont\dbfam= \sevendb
\def\db {\fam\dbfam\tendb}
\font\tenrsfs=rsfs10 
\font\sevenrsfs=rsfs7
\newfam\scrfam 
\textfont\scrfam = \tenrsfs
\scriptfont\scrfam= \sevenrsfs

\def\N{{\cal N}}

\def\PP{{\cal P}}

\def\f{{\bi f}}

\def\F{{\bi F}}

\def\Q{{\db Q}}

\def\C{{\db C}}

\def\F{{\db F}}

\def\Gal{\mathop{\rm Gal}\nolimits}

\font\tendb=msbm10 
\font\sevendb=msbm7
\newfam\dbfam 
\textfont\dbfam = \tendb 
\scriptfont\dbfam= \sevendb
\def\db {\fam\dbfam\tendb}
\font\tenrsfs=rsfs10 
\font\sevenrsfs=rsfs7
\newfam\scrfam 
\textfont\scrfam = \tenrsfs
\scriptfont\scrfam= \sevenrsfs

\def\f{{\bi f}}

\def\Q{{\db Q}}

\def\N{{\db N}}

\def\GL{{\db GL}}

\def\F{{\db F}}

\def\Gal{\mathop{\rm Gal}\nolimits}
\def\mod{\hbox{\rm \ mod.\ }}

\font\teufm=eufm10
\font\seufm=eufm10 at 7pt
\font\sseufm=eufm10 at 6pt
\newfam\fameufm
\textfont\fameufm=\teufm
\scriptfont\fameufm=\seufm
\scriptscriptfont\fameufm=\sseufm
\def\goth{\fam\fameufm\teufm}

\def\wp{{\goth p}}

\def\qq{{\goth q}}

\def\PP{{\goth P}}
\def\m{{\goth m}}

\def\f{{\goth f}}

\def\NN{{\goth N}}

\def\S{{\cal S}}

\def\Gal{\mathop{\rm Gal}\nolimits}

\def\NN{\mathop{\rm N}\nolimits}

\def\mod{\hbox{\rm \ mod.\ }}

\bigskip\bigskip

\bigskip\bigskip
\def\nn{{\goth N}}

\def\N{{\cal N}}

\def\E{{\cal E}}

\def\og{\leavevmode\raise.3ex\hbox{$\scriptscriptstyle\langle\!\langle\,$}}
\def\fg{\leavevmode\raise.3ex\hbox{$\scriptscriptstyle\,\rangle\!\rangle$}}

\centerline {\titre{Sur le th\'eor\`eme de Fermat sur $\Q\bigl(\sqrt{5}\bigr)$}}
\bigskip
\smallskip
\centerline {Alain Kraus}
\bigskip 
\smallskip
{\bf{Abstract.}} Let $p$ be an odd prime number. Using modular arguments, we give 
an easy testable condition which allows often  to prove  Fermat's Last Theorem over  the quadratic field $\Q\bigl(\sqrt{5}\bigr)$ for the exponent $p$. 
It is related to the Wendt's resultant of the polynomials $X^n-1$ and $(X+1)^n-1$. 
We  deduce Fermat's Last Theorem over this field in case one has
 $5\leq p<10^7$, and we  obtain  analogous results on Sophie Germain type criteria.
\bigskip

\medskip

\bigskip
\medskip

\centerline{INTRODUCTION}
\bigskip
Soient $p$ un nombre premier impair et $K$ un corps de nombres.  On dit que le  th\'eor\`eme de Fermat est vrai sur $K$ pour 
l'exposant $p$, s'il n'existe pas  d'\'el\'ements $a,b,c$ dans  $K$ tels que  l'on ait
$$a^p+b^p+c^p=0\quad \hbox{et}\quad  abc\neq 0.$$
En 1994, A. Wiles a d\'emontr\'e  qu'il en est ainsi pour  $K=\Q$ ([18]). Dix ans plus tard, F. Jarvis et P. Meekin ont \'etendu ce r\'esultat au corps $\Q\bigl(\sqrt{2}\bigr)$ ([8]). 
R\'ecemment, N. Freitas et S. Siksek ont effectu\'e des progr\`es importants sur cette \'equation  si $K$ est un corps totalement r\'eel. Ils obtiennent  dans [3] un crit\`ere pour que le th\'eor\`eme de Fermat soit vrai sur $K$ asymptotiquement i.e.  si $p$ est plus grand qu'une constante ne d\'ependant  que de $K$. Cela leur permet de d\'emontrer le th\'eor\`eme de Fermat asymptotique pour  une proportion de cinq sixi\`eme des corps quadratiques r\'eels. Par ailleurs, ils d\'emontrent dans [4] le th\'eor\`eme de Fermat sur les corps $\Q\bigl(\sqrt{d}\bigr)$ avec   $3\leq d\leq 23$  et $d\neq 5,17$, pour   tout $p\geq 5$. 
Ils parviennent  aussi \`a conclure  sur $\Q\bigl(\sqrt{17}\bigr)$ pour une infinit\'e explicite de nombres premiers.
\smallskip
On s'int\'eresse ici \`a l'\'etude de l'\'equation de Fermat sur  $\Q\bigl(\sqrt{5}\bigr)$. Les r\'esultats obtenus dans [3]
ne permettent pas,  a priori, de d\'emontrer le th\'eor\`eme de Fermat asymptotique sur   ce corps.
Il semble que les r\'esultats d\'ej\`a connus soient les suivants.  L'\'equation de Fermat d'exposant $3$ est une courbe elliptique ayant une infinit\'e de points rationnels 
 sur  $\Q\bigl(\sqrt{5}\bigr)$. On a par exemple l'\'egalit\'e
 $$\bigl(9+\sqrt{5}\bigr)^3+\bigl(9-\sqrt{5}\bigr)^3=12^3.$$
 Par ailleurs, le th\'eor\`eme  de Fermat   est vrai sur $\Q\bigl(\sqrt{5}\bigr)$ pour $p=5,7,11$ ([6], th. 5.1]),  et  pour les nombres premiers $p$
congrus \`a $\pm 1$ modulo $5$  qui ne divisent pas le nombre de classes du corps $\Q\big(\sqrt{5},\mu_p\bigr)$, o\`u $\mu_p$ est le groupe des racines $p$-i\`emes de l'unit\'e ([7], th. 4).  
\smallskip
\'Etant donn\'e un nombre premier $p\geq 11$, on obtient ici, par des m\'ethodes modulaires, un crit\`ere simple permettant souvent en pratique de d\'emontrer le th\'eor\`eme de Fermat  sur $\Q\bigl(\sqrt{5}\bigr)$ pour l'exposant $p$.  Il peut   sans doute s'adapter \`a d'autres corps quadratiques. Je n'ai pas \'et\'e plus loin dans cette direction.
\bigskip

{\bf{ 1. \'Enonc\'e des r\'esultats}}
\medskip
Pour tout entier $n\geq 1$, notons $W_n$ le r\'esultant des polyn\^omes 
$X^n-1$ et $(X+1)^n-1.$
 La lettre $p$ d\'esigne un nombre premier $\geq 5$. Posons d\'esormais
$K=\Q\bigl(\sqrt{5}\bigr).$
\bigskip

\proclaim Th\'eor\`eme. Supposons qu'il existe un entier $n\geq 1$  satisfaisant les conditions suivantes~:
 \smallskip
\vskip0pt\noindent
1) on a $n<p-2$ et $n\equiv 2 \mod 4$.
\smallskip
\vskip0pt\noindent
2) L'entier $np+1$ est un nombre premier, congru \`a $\pm 1$ modulo $5$, ne divisant pas $W_n$.
\smallskip
\vskip0pt\noindent
Alors, le th\'eor\`eme de Fermat est vrai sur  $K$ pour l'exposant $p$.
\bigskip

Des exp\'erimentations num\'eriques effectu\'ees \`a l'aide du logiciel de calculs Pari [1],  rendent plausible le fait qu'il existe toujours un tel entier $n$ si  $p$ n'est pas 
 dans l'ensemble
$$\Big\lbrace 7,11,23,53,59,67,79,83,127\Big\rbrace.$$
Pour autant,   Dickson  a d\'emontr\'e en 1909 que 
tout nombre premier de la forme $np+1$ assez grand, par exemple plus grand que $p^4$, 
 divise $W_n$ ([14], p. 137).  La  question de savoir  s'il existe  un entier $n$ tel que $np+1$ soit un nombre premier   ne divisant  pas $W_n$ est  ouverte.
 \bigskip

\proclaim Corollaire 1. Si $p$ est plus petit que $10^7$, le th\'eor\`eme de Fermat est vrai sur $K$ pour l'exposant $p$.
\bigskip

\`A titre indicatif, la liste des couples $(p,n)$ pour $p<50$, 
 avec les plus petits entiers $n$ pour lesquels le crit\`ere fonctionne est la suivante :
 $$(5,2), (13,10), (17,14), (19,10), (29,2), (31,10), (37,34), (41,38), (43,10), (47,14).$$
Ce crit\`ere  permet  aussi d'\'etablir le th\'eor\`eme de Fermat sur  $K$ pour des $\og$grands$\fg$ exposants. Par exemple, pour $p=10^{100}+267$, qui est le plus petit nombre premier plus grand que $10^{100}$, on peut conclure  avec $n=754$.
 \bigskip
 
\proclaim Corollaire 2. Supposons que l'une des conditions suivantes soit satisfaite :
\smallskip
\vskip0pt\noindent
1) on a $p\equiv 4 \mod 5$ et  $2p+1$ est premier.
\smallskip
\vskip0pt\noindent
2) L'entier  $10p+1$ est premier. 
\smallskip
\vskip0pt\noindent
Alors, le th\'eor\`eme de Fermat est vrai sur  $K$ pour l'exposant $p$.
\bigskip

{\bf 2. {Principe de d\'emonstration}}
\medskip
Il repose sur la m\'ethode modulaire, qui est expos\'ee en d\'etail dans [3] lorsque le corps de base est  totalement r\'eel.  Supposons dans toute la suite $p\geq 11$, cela n'est pas restrictif.
\smallskip
On proc\`ede par l'absurde en supposant 
qu'il existe un triplet $(a,b,c)$ d'\'el\'ements de $K$  tel que $a^p+b^p+c^p=0$ et $abc\neq 0$. On peut supposer que $a,b,c$ sont dans l'anneau d'entiers $O_K$ de $K$. 
Parce que  $O_K$ est principal, on peut de plus supposer  $a$, $b$, $c$  premiers entre eux.
Soit $E_0$ la courbe elliptique d\'efinie sur $K$ d'\'equation
$$y^2=x(x-a^p)(x+b^p).$$
Elle   est semi-stable en dehors de  l'id\'eal premier $\PP_2$ de $O_K$ au-dessus de $2$. De plus, pour tout id\'eal premier $\qq$ de $O_K$, distinct de $\PP_2$, 
l'exposant en $\qq$ de l'id\'eal discriminant minimal de $E_0$  est divisible par $p$. 
En normalisant $(a,b,c)$ convenablement, on d\'emontre que le conducteur local de $E_0$ en $\PP_2$  est
$$\PP_2^r\quad \hbox{avec}\quad r\in\Big\lbrace 1,2,3\Big\rbrace.$$
Par ailleurs, $E_0$ est modulaire ([5]). 
Soit $\overline \Q$ la fermeture alg\'ebrique de $\Q$ dans $\C$. Notons 
$$\rho_{E_0,p} : \Gal\bigl(\overline \Q/K\bigr)\to \GL_2(\F_p)$$
 la repr\'esentation donnant, moyennant le choix d'une base, l'action du groupe de Galois  
$\Gal\bigl(\overline \Q/K\bigr)$ sur le 
groupe des points de $p$-torsion de $E_0$. 
En utilisant la th\'eorie du corps de classes, on v\'erifie qu'elle est irr\'eductible.  

Il r\'esulte alors du th\'eor\`eme d'abaissement du niveau ([3], th. 7), que $\rho_{E_0,p}$ $\og$provient$\fg$, en un sens \`a pr\'eciser,
 d'une newform  modulaire  parabolique de Hilbert de poids $(2,2)$ et de niveau 
$\PP_2^r$.
Or il  n'existe pas de telles formes si $r=1$ et $r=2$ (cf. [12]).  Il en existe une seule si $r=3$ ({\it{loc. cit.}}), qui corres\-pond \`a une courbe elliptique $E$ d\'efinie sur $K$, de conducteur $\PP_2^3$, unique \`a $K$-isog\'enie pr\`es ([17]).  On en d\'eduit  que les modules galoisiens des points de $p$-torsion de $E_0$ et  $E$ sont isomorphes. Sous les hypoth\`eses du th\'eor\`eme, cela permet d'en d\'eduire  une contradiction \`a  l'existence  de $(a,b,c)$.
\bigskip

{\bf{3. L'anneau $O_K/\PP_2^2$}}
\medskip
Posons 
 $$u={1+\sqrt{5}\over 2}.$$
On a  $u^2=1+u$. L'ensemble 
 $$\Big\lbrace 0,1,u,u^2\Big\rbrace$$
  est  un syst\`eme de repr\'esentants du  corps r\'esiduel $O_K/\PP_2$. L'anneau $O_K/\PP_2^2$ est  local de cardinal 16 et on a 
 $$O_K/\PP_2^2=\Big\lbrace x+2y+\PP_2^2\ \big|\ x,y=0,1,u,u^2\Big\rbrace.$$
Le groupe  de ses  \'el\'ements inversibles est   
  $$\Big\lbrace x+2y+\PP_2^2\ \big|\ x=1,u,u^2,\ y=0,1,u,u^2   \Big\rbrace.$$
 \smallskip

{\bf{ 4. Normalisation des solutions}}
\medskip
Soit  $(a,b,c)$ un triplet d'\'el\'ements de $O_K$, tous non nuls et premiers entre eux, tel que $a^p+b^p+c^p=0$. 
\smallskip

\proclaim Lemme 1. Quitte \`a  permuter  $a,b,c$ entre eux et \`a multiplier $(a,b,c)$ par une unit\'e convenable de $O_K$,  l'une des deux conditions 
suivantes est satisfaite~: 
\smallskip
\vskip0pt\noindent
1) $2$ divise $a$, sans diviser $bc$, et  $b$ est un carr\'e modulo $4$.
\smallskip
\vskip0pt\noindent
2) $2$ ne divise pas $abc$ et $(a^p,b^p)$ est congru modulo $4$ \`a l'un des couples  
$$(u,1+2u),\quad (3,u^2),\quad (1,u), \quad (1,u^2+2u).$$

D\'emonstration :  1) Supposons que  $2$ divise $abc$.  Parce que $a,b,c$ sont premiers entre eux,  on peut supposer que $2$ divise $a$ sans diviser $bc$.
Il s'agit alors de v\'erifier qu'il existe une unit\'e $\xi$ de $O_K$ telle que $\xi b$ soit un carr\'e modulo $4$. On a $b\equiv 1,u,u^2\mod 2$. Quitte \`a multiplier  $(a,b,c)$ par $-1$, on peut supposer que l'on a
$$b\equiv 1,u+2u^2,u^2\mod 4\quad \hbox{ou bien}\quad b\equiv 1+2u,u,u^2+2 \mod 4.$$
On  a $1+2u=u^3$ et $u(u^2+2)=1+4u$. On en d\'eduit respectivement que $\xi=1$ ou bien que $\xi=u$ convient, d'o\`u la condition 1.
\smallskip
2) Supposons que $2$ ne divise pas $abc$. On commence par normaliser $(a,b,c)$ de sorte que l'on ait
$$a\equiv 1 \mod 4\quad \hbox{et}\quad b^p\equiv u^2 \mod 2.\leqno(1)$$
Pour cela, on remarque que l'un des entiers $a,b,c$ est congru \`a $1$ modulo $2$. Supposons que ce soit $a$. Quitte \`a multiplier $(a,b,c)$ par $-1$, on peut  alors supposer (cf. l'alin\'ea 1)
$$a\equiv 1\mod 4\quad \hbox{ou}\quad a\equiv 1+2u\mod 4.$$
Si on a $a\equiv 1+2u\mod 4$, en multipliant $(a,b,c)$ par $1+2u$, on se ram\`ene au cas o\`u $a\equiv 1 \mod 4$.  Quitte \`a \'echanger $b$ et $c$, la condition (1) est  alors satisfaite. 
\smallskip
On peut donc supposer au d\'epart que  l'on a 
$$a^p\equiv 1 \mod 4 \quad \hbox{et}\quad b^p\equiv u^2,3u^2, u^2+2, u^2+2u\mod 4.\leqno(2)$$
Par ailleurs, $u$ est   d'ordre $6$ modulo $4$, d'o\`u 
$$u^p\equiv u \mod 4 \ \hbox{si}\ p\equiv 1 \mod 6\quad \hbox{et}\quad u^p\equiv u^{-1}\equiv u^2+2 \mod 4 \ \hbox{si}\ p\equiv -1 \mod 6.\leqno(3)$$

2.1) Si on a $b^p\equiv u^2 \mod 4$, d'apr\`es (2) et (3), le couple $\bigl((u a)^p,(ub)^p\bigr)$ est donc congru modulo $4$ \`a 
$$(u,1+2u)\quad \hbox{ou} \quad (u^2+2,u).$$
Par ailleurs, on a $u^2(u^2+2,u)\equiv (u,1+2u)\mod 4$. Vu que $u^2$ est congru \`a $u^{2p}$ ou \`a $u^{4p}$ modulo $4$, 
on peut ainsi supposer   $(a^p,b^p)\equiv (u,1+2u)\mod 4$.

2.2) Si on a  $b^p\equiv 3u^2 \mod 4$, en multipliant  $(a,b,c)$ par  $-1$, on se ram\`ene au cas o\`u $(a^p,b^p)\equiv (3,u^2)\mod 4$.

2.3) Si on a  $b^p\equiv u^2+2\mod 4$, alors $\bigl((u b)^p,(ua)^p\bigr)$ est  congru modulo $4$ \`a 
$$(1,u)\quad \hbox{ou}\quad (u+2u^2,u^2+2).$$
On a $u^2(u+2u^2,u^2+2)\equiv (1,u)\mod 4$, et  dans ce cas on peut  supposer $(a^p,b^p)$ congru \`a $(1,u)$ modulo $4$.
 
2.4) Si on a $b^p\equiv u^2+2u \mod 4$, on obtient le dernier couple de la liste annonc\'ee, d'o\`u le r\'esultat.

\bigskip 
{\bf{5. La courbe elliptique $E_0$}}
\medskip
Soit  $(a,b,c)$ un triplet d'\'el\'ements de $O_K$, tous non nuls et premiers entre eux, tel que $a^p+b^p+c^p=0$, normalis\'e de sorte que 
l'une des deux conditions du lemme 1  soit satisfaite.
Pour toute la suite,  posons
$$A=a^p,\quad B=b^p,\quad C=c^p,$$
et notons $E_0$  la cubique de Weierstrass d'\'equation 
$$y^2=x(x-A)(x+B).\leqno(4)$$
Les invariants standard $c_4$, $c_6$ et $\Delta$ associ\'es \`a cette \'equation sont
$$c_4=16\bigl(A^2+AB+B^2\bigr),\quad c_6=-32\bigl(A-B\bigr)\bigl(B-C\bigr)\bigl(C-A\bigr),\quad \Delta=16(ABC)^2.$$
En particulier, $E_0$ est une courbe elliptique d\'efinie sur $K$.   Pour tout id\'eal premier non nul de $O_K$, notons  $v_{\qq}$ la valuation sur $K$ qui lui est associ\'ee. 
\bigskip

\proclaim Lemme 2. Soit $\qq$ un id\'eal premier de $O_K$ distinct de $\PP_2$.
\smallskip
\vskip0pt\noindent
 1) L'\'equation $(4)$ est   minimale en $\qq$.
\smallskip
\vskip0pt\noindent
 2) Si $v_{\qq}(abc)=0$, $E_0$ a bonne r\'eduction en $\qq$. 
\smallskip
\vskip0pt\noindent
 3) Si $v_{\qq}(abc)\geq 1$, $E_0$ a  r\'eduction de type multiplicatif en $\qq$ et $p$ divise $v_{\qq}(\Delta)$.

D\'emonstration :  Si   $v_{\qq}(abc)=0$, alors par d\'efinition $E_0$ a bonne r\'eduction en $\qq$. Si  $v_{\qq}(abc)\geq 1$,  $a,b,c$ \'etant  premiers entre eux, on a $v_{\qq}(c_4)=0$, d'o\`u  le lemme.
\bigskip

Notons $\N_0$ le  conducteur de $E_0$.
\bigskip

\proclaim Lemme 3. 1)  Si $2$ divise $abc$, $E_0$ a   r\'eduction de type multiplicatif en $\PP_2$.
\smallskip
\vskip0pt\noindent
2) Supposons que  $2$ ne divise pas $abc$. Alors $E_0$ a  r\'eduction de type additif en $\PP_2$. On a $v_{\PP_2}(\N_0)\in \big\lbrace 2,3\big\rbrace$ et $E_0$ a potentiellement bonne r\'eduction en $\PP_2$.

D\'emonstration : 1) Supposons que $2$ divise $abc$. On a 
$v_{\PP_2}(a)\geq 1$ et $v_{\PP_2}(bc)=0$. 
D'apr\`es l'\'egalit\'e $c_4^3-c_6^2=1728\Delta$, on a donc 
$$v_{\PP_2}(c_4)=4,\quad v_{\PP_2}(c_6)=6,\quad v_{\PP_2}(\Delta)\geq 26.\leqno(5)$$
V\'erifions alors que l'\'equation (4) n'est pas minimale en $\PP_2$, ce qui d'apr\`es (5) \'etablira la premi\`ere assertion.
Suivant la terminologie utilis\'ee dans [13], le type de r\'eduction de $E_0$ en $\PP_2$ correspond au cas 7 de Tate ou bien l'\'equation (4) n'est pas minimale en $\PP_2$ (tableau IV, p. 129 ; dans la colonne $\og$Equation non minimale$\fg$, le triplet de valuations $(4,6,12)$ doit \^etre remplac\'e par $(4,6\geq 12)$). Les invariants standard $b_2, b_4, b_6, b_8$ associ\'es \`a l'\'equation (4) sont
$$b_2=4\bigl(B-A\bigr),\quad b_4=-2AB,\quad b_6=0,\quad b_8=-(AB)^2.$$
Avec les notations de la proposition 4 de [13], l'entier  $r=0$ satisfait la condition (a). Par ailleurs, $b$ est un carr\'e modulo $4$ (lemme 1), donc $B$ aussi. On en d\'eduit que l'on est dans un cas de Tate $\geq 8$, ce qui implique l'assertion.
\smallskip
2) Supposons que $2$ ne divise pas $abc$. 
 \smallskip
2.1) Si  $(A,B)$ est congru \`a $(u,1+2u)$ modulo $4$, on v\'erifie que l'on~a 
 $$v_{\PP_2}(c_4)=6, \quad v_{\PP_2}(c_6)=5,\quad   v_{\PP_2}(\Delta)=4.$$
Avec les notations de la proposition 1 de [13], 
en choisissant    $(r,t)=(u^2,1)$, on constate que le type de r\'eduction de $E_0$ en $\PP_2$ correspond au cas 5 de Tate, 
d'o\`u $v_{\PP_2}(\N_0)=2.$
  \smallskip
2.2) Pour  les  autres congruences   de $(A,B)$ modulo $4$ figurant dans l'\'enonc\'e du lemme 1,
 on v\'erifie que l'on a 
 $$v_{\PP_2}(c_4)=5, \quad v_{\PP_2}(c_6)=5,\quad   v_{\PP_2}(\Delta)=4.$$
 Si $(A,B)$ est congru \`a  $(3,u^2)$ ou  $(1,u^2+2u)$ modulo $4$, on constate que le type de r\'eduction de $E_0$ en $\PP_2$ correspond au cas 4 de Tate ([13], prop. 1, avec $r=u$ et $t=1$), d'o\`u $v_{\PP_2}(\N_0)=3.$
\smallskip
\vskip0pt\noindent
Si on a $(A,B)\equiv (1,u)\mod 4$, on aboutit \`a la m\^eme conclusion 
avec  $(r,t)=(u^2,1)$.
 \smallskip
\vskip0pt\noindent
 Dans chacun des cas envisag\'es, $E_0$ a un invariant modulaire entier en $\PP_2$, donc $E_0$ a potentiellement bonne r\'eduction en $\PP_2$, d'o\`u le r\'esultat.
 \bigskip

{\bf{Remarque.}} Si $2$ ne divise pas $abc$, on ne peut pas toujours, a priori, normaliser $(a,b,c)$ de sorte que l'on ait  $v_{\PP_2}(\N_{0})=2$. En effet, 
supposons par exemple  
$$A\equiv 1 \mod 4\quad \hbox{et}\quad   B\equiv  u^2+2u \mod 4.$$
On a alors $C\equiv u+2 \mod 4$ et les \'egalit\'es
$$v_{\PP_2}\bigl(A^2+AB+B^2\bigr)=v_{\PP_2}\bigl(A^2+AC+C^2\bigr)=v_{\PP_2}\bigl(B^2+BC+C^2\bigr)=1.$$
Par suite, la valuation en $\PP_2$ de l'invariant $c_4$ d'un mod\`ele minimal de toute courbe elliptique d\'eduite de $E_0$ apr\`es normalisation de $(a,b,c)$, vaut $5$. L'exposant en $\PP_2$ de son conducteur est donc distinct de $2$. 
\bigskip

 {\bf{ 6. La repr\'esentation $\rho_{E_0,p}$}}
 \medskip
L'\'enonc\'e qui suit est en particulier d\'emontr\'e dans [4], tout au moins si $p\geq 17$ et avec la norma\-lisation  adopt\'ee dans {\it{loc. cit.}}  pour $(a,b,c)$. On se limitera ici \`a en rappeler bri\`evement les arguments et \`a en donner  une d\'emonstration  plus simple dans notre situation. Notons d\'esormais $S_p$ l'ensemble des id\'eaux premiers de $O_K$ au-dessus de $p$. Rappelons que l'on a suppos\'e $p\geq 11$. 
\bigskip

\proclaim Proposition.  La repr\'esentation $\rho_{E_0,p}$ est irr\'eductible.

 D\'emonstration :  Supposons que $\rho_{E_0,p}$ soit  r\'eductible. Il existe deux caract\`eres $\varphi, \varphi'$ de  $\Gal\bigl(\overline \Q/K\bigr)$ \`a valeurs dans $\F_p^*$ tels que 
 $\rho_{E_0,p}$ soit repr\'esentable sous la forme
 $$\pmatrix{\varphi & *\cr 0&\varphi'\cr}.$$
 Les caract\`eres  $\varphi$ et  $\varphi'$ sont non ramifi\'es en dehors de $S_p\cup\big\lbrace \PP_2\big\rbrace$ ; 
de plus, pour tout $\wp\in S_p$, la courbe elliptique $E_0$ est semi-stable en $\wp$, donc $\varphi$ ou  $\varphi'$ est non ramifi\'e en $\wp$ (cf. [15], p. 274-277]).
\smallskip
 1) Supposons que $2$ divise $abc$. Dans ce cas, $E_0$ est semi-stable et  $\varphi, \varphi'$ sont donc non ramifi\'es en dehors de $S_p$ ([11], prop. 3). Si $p$ est inerte dans $K$, $\varphi$ ou $\varphi'$ est ainsi  partout non ramifi\'e aux places finies de $K$. Le nombre de classes restreint de $K$ vaut $1$. Par suite,  $\varphi$ ou $\varphi'$ est trivial. On peut supposer que c'est $\varphi$, auquel cas $E_0$ poss\`ede un point d'ordre $p$ rationnel sur $K$. Tous ses points d'ordre 2 sont dans $K$, en particulier, $E_0$ a donc un point d'ordre $2p$ rationnel sur $K$. On obtient  une contradiction car    $p\geq 11$ ([9], th. 1). Supposons  $p$  d\'ecompos\'e dans $K$. Posons
$pO_K=\wp \wp'$. D'apr\`es  ce qui pr\'ec\`ede, 
on peut supposer  $\varphi$ ramifi\'e en $\wp$ et non ramifi\'e en $\wp'$, auquel cas la restriction de $\varphi$ \`a un sous-groupe d'inertie en $\wp$  est le caract\`ere cyclotomique. Cela implique $u^2\equiv 1 \mod \wp$ et de nouveau une contradiction ([11], th. Appendice A), d'o\`u l'assertion dans ce cas.
\smallskip
 2) Supposons que $2$ ne divise pas $abc$. Compte tenu du lemme 3, les images de $\varphi$ et $\varphi'$ restreintes \`a un sous-groupe d'inertie en $\PP_2$ sont d'ordre divisible par $3$ (cf. [15], prop. 23 et [10], th. 3). 
  
2.1) Supposons $\varphi$ ou $\varphi'$ non ramifi\'e en dehors de $\PP_2$, par exemple $\varphi$,  ce n'est pas restrictif. Soit $L$ le corps laiss\'e fixe par $\varphi$. C'est une extension cyclique de $K$.  Son conducteur est une puissance de $\PP_2$. Il existe donc $r\geq 1$ tel que $L$ soit contenu dans le corps de rayon $K^{\m}$ avec  $\m=\infty_1\infty_2\PP_2^r$, o\`u $\infty_1$ et $\infty_2$ sont les deux places \`a l'infini de $K$.   Par suite, $3$   divise le degr\'e de $K^{\m}$ sur $K$. V\'erifions  en fait que ce degr\'e  est une puissance de $2$, ce qui entra\^\i nera la contradiction cherch\'ee.
Notons $U_K$ le groupe des unit\'es de $O_K$ et $U_{\m,1}$ le sous-groupe de $U_K$ form\'e des unit\'es congrues \`a $1$ modulo $\m$. Le degr\'e de $K^{\m}$ sur $K$ est (cf. [2], cor. 3.2.4)
$${4^{r}\times 3\over (U_K:U_{\m,1})}.$$
 Il s'agit donc  de montrer que $3$ divise l'indice $(U_K:U_{\m,1})$. On v\'erifie pour cela que  $u^2$ est d'ordre $3$ modulo $\PP_2$, d'o\`u l'on d\'eduit que $3$ divise l'ordre de $u^2$ modulo $\PP_2^r$. Parce que $u^2$ est totalement positive, cela entra\^\i ne notre assertion.

2.2) Supposons  $\varphi$  et $\varphi'$ ramifi\'es en dehors de $\PP_2$. Dans ce cas,  $p$ est n\'ecessairement d\'ecompos\'e dans $K$. Posons $pO_K=\wp \wp'$.  On peut  supposer $\varphi$ ramifi\'e en $\wp$ et non ramifi\'e en $\wp'$. Il en r\'esulte que 
 $\varphi^{6}$ est non ramifi\'e en dehors de $\wp$ et que sa restriction \`a un sous-groupe d'inertie en $\wp$ est la puissance sixi\`eme du caract\`ere cyclotomique ([4], p. 11). On en d\'eduit que l'on a  $u^{12}\equiv 1 \mod \wp$ ({\it{loc. cit.}}),  ce qui conduit \`a  une contradiction, d'o\`u  le r\'esultat.
\bigskip

 {\bf{ 7. D\'emonstrations des r\'esultats annonc\'es}}
 \medskip
Pour tout id\'eal non nul $\nn$ de $O_K$, notons $\NN(\nn)$ le cardinal de $O_K/\nn$. \'Etant donn\'ee une courbe elliptique $\E$ d\'efinie sur $K$, on d\'esignera par  $$L(\E,s)=\sum_{\nn} {a_{\nn}(\E)\over \NN(\nn)^{s}}$$
 sa fonction $L$ de Hasse-Weil. Pour tout id\'eal premier  $\qq$ de $O_K$ en lequel $\E$ a bonne r\'eduction,  
 $$\NN(\qq)+1-a_{\qq}(\E)$$
est le nombre de points  rationnels sur $O_K/\qq$ de la courbe elliptique d\'eduite de $\E$ par r\'eduction.  Rappelons par ailleurs que $\E$ est modulaire ([5]).
 \smallskip

Supposons  
que le th\'eor\`eme de Fermat soit faux sur $K$ pour l'exposant $p$. Il existe alors $a,b,c$ dans $O_K$, 
tous non nuls et premiers entre eux, satisfaisant l'une des conditions du lemme 1, tels que $a^p+b^p+c^p=0$.  Soit $E_0$ la courbe elliptique  associ\'ee \`a $(a,b,c)$ comme dans le paragraphe 5. 
\smallskip
Pour $r\in \big\lbrace 1,2,3\big\rbrace$, notons 
$\S_2(\PP_2^r)$ l'ensemble des newforms modulaires paraboliques de Hilbert de poids parall\`ele $(2,2)$ et de niveau $\PP_2^r$. C'est un $\C$-espace vectoriel de dimension finie. 
Compte tenu des lemmes 2 et 3, et de la proposition,  on d\'eduit  du  th\'eor\`eme d'abaissement du niveau ([3], th. 7) l'existence d'une newform
$\f\in \S_2(\PP_2^r)$
et d'une place $\PP$ de $\overline \Q$ de ca\-ract\'eristique r\'esiduelle $p$ 
telles que, pour tout id\'eal premier $\qq$ de $O_K$  qui n'est pas dans $S_p\cup \big\lbrace \PP_2\big\rbrace$, 
en notant $a_{\qq}(\f)$ le coefficient de Fourier de $\f$ en $\qq$, les conditions suivantes soient satisfaites :
$$a_{\qq}(\f)\equiv a_{\qq}(E_0) \mod \PP \quad \hbox{si}\quad E_0 \ \hbox{a bonne r\'eduction en}\ \qq,$$
$$a_{\qq}(\f)\equiv \pm \bigl(\NN(\qq)+1\bigr) \mod \PP \quad \hbox{si}\quad E_0 \ \hbox{a  r\'eduction de type multiplicatif en}\ \qq.$$
\bigskip

 {\bf{7.1. Le  th\'eor\`eme}}
 \medskip
Les $\C$-espaces vectoriels $\S_2(\PP_2)$ et $\S_2(\PP_2^2)$ sont nuls (cf. [12]). On a donc $r=3$.
Le $\C$-espace vectoriel $\S_2(\PP_2^3)$ est  de dimension 1 ({\it{loc. cit.}}).
D'apr\`es les tables de [17], il existe une courbe elliptique $E$ d\'efinie sur $K$, unique \`a $K$-isog\'enie pr\`es, de conducteur $\PP_2^3$. Elle est modulaire. Par suite, les fonctions $L$ de $E$ et de la forme modulaire de Hilbert normalis\'ee de $\S_2(\PP_2^3)$ sont \'egales.
Ainsi, pour tout id\'eal premier $\qq$ de $O_K$,  qui n'est pas dans $S_p\cup \big\lbrace \PP_2\big\rbrace$, on  a
$$a_{\qq}(E)\equiv a_{\qq}(E_0) \mod p \quad \hbox{si}\quad E_0 \ \hbox{a bonne r\'eduction en}\ \qq,\leqno(6)$$
$$a_{\qq}(E)\equiv \pm \bigl(\NN (\qq)+1\bigr) \mod p \quad \hbox{si}\quad E_0 \ \hbox{a  r\'eduction de type multiplicatif en}\ \qq.\leqno(7)$$
Consid\'erons un entier $n$ satisfaisant les deux conditions de l'\'enonc\'e du th\'eor\`eme et posons $q=np+1$.  Soit $\qq$ un id\'eal premier de $O_K$ au-dessus de $q$. Il n'est pas dans $S_p\cup \big\lbrace \PP_2\big\rbrace$.
\bigskip

\proclaim Lemme 4.  La courbe elliptique $E_0$ a  r\'eduction de type multiplicatif en $\qq$.

D\'emonstration :  Supposons le contraire. Dans ce cas, $E_0$ a bonne r\'eduction en $\qq$ et on a $v_{\qq}(abc)=0$. Parce que $q$ est d\'ecompos\'e dans $K$, le corps $O_K/\qq$ est de cardinal $q$. Il en r\'esulte que  $a^p$, $b^p$ et  $c^p$ sont  des racines $n$-i\`emes de l'unit\'e modulo $\qq$. Posons 
$$\alpha=\Bigl({a\over c}\Bigr)^p \mod \qq.$$
On a $\alpha^n\equiv 1 \mod \qq$. D'apr\`es l'\'egalit\'e $a^p+b^p+c^p=0$, on a $\alpha+1\equiv -\bigl({b\over c}\bigr)^p \mod \qq$, d'o\`u
$(\alpha+1)^n\equiv 1 \mod \qq$ ($n$ est pair). 
Les polyn\^omes $X^n-1$ et $(X+1)^n-1$ ont donc une racine commune modulo $\qq$, ce qui contredit le fait que $q$ ne divise pas $W_n$, d'o\`u le lemme.
\medskip
 
D'apr\`es la condition (7) et le lemme 4, on a donc
$$a_{\qq}(E)\equiv \pm 2 \mod p.$$
Par ailleurs,  $E$ a bonne r\'eduction en $\qq$ et  a tous ses points d'ordre 2 rationnels sur $K$ ([17]). On a donc  (cf. [16], p. 192, prop. 3.1)
$$a_{\qq}(E)\equiv q+1 \mod 4.\leqno(8)$$
 En particulier, $a_{\qq}(E)$ est pair et on a ainsi
$$a_{\qq}(E)\equiv \pm 2 \mod 2p.$$
Puisque $q$ est d\'ecompos\'e dans $K$,  on a d'apr\`es les bornes de Weil
$$|a_{\qq}(E)|\leq 2\sqrt{q}.$$
Supposons $a_{\qq}(E)\neq \pm 2$. L'entier $2p$ divise $a_{\qq}(E)\pm 2$, d'o\`u alors l'in\'egalit\'e
$$2p\leq 2\sqrt{q}+2\quad \hbox{i.e.}\quad p\leq \sqrt{q}+1.$$
Cela implique
$$p\leq n+2,$$
ce qui conduit \`a une contradiction.  On a donc 
$$a_{\qq}(E)=\pm 2,$$
et $q+1-a_{\qq}(E)$ vaut $q+3$ ou $q-1$. D'apr\`es  la condition (8), on en d\'eduit que l'on a 
$$q\equiv 1 \mod 4.$$  Or  par hypoth\`ese, on a  $n\equiv 2 \mod 4$, d'o\`u 
$q\equiv 2p+1\equiv 3 \mod 4$. On obtient ainsi une contradiction \`a l'existence de $(a,b,c)$. Cela termine la d\'emonstration du th\'eor\`eme. 
\bigskip

{\bf{ 7.2. Les corollaires 1 et 2}}
 \medskip
1) On v\'erifie \`a l'aide du logiciel de calculs Pari ([1]), que si $p$ est plus petit que $10^7$ et est distinct de 
$$11,23,53,59,67,79,83,127,$$
alors 
$p$ satisfait les deux conditions de l'\'enonc\'e du th\'eor\`eme, d'o\`u le r\'esultat dans ce cas (le temps d'ex\'ecution qui m'a \'et\'e n\'ecessaire sur un MacBook Pro est de 4h30mn).
\smallskip

En ce qui concerne ces huit nombres premiers, on proc\`ede  pour conclure de fa\c con analogue. On d\'etermine un entier $n$ convenable de sorte que $q=np+1$ soit un nombre premier d\'ecompos\'e dans $K$ ne divisant pas $W_n$,  \`a  ceci pr\`es  que la condition 1  n'est pas satisfaite.
\vskip0pt\noindent
 Plus pr\'ecis\'ement, supposons $p=11$.  On choisit  $n=8$. L'entier $q=89$ est un nombre premier d\'ecompos\'e dans $K$. Si $\qq$ est un id\'eal premier de $O_K$ au-dessus de $89$, on constate que l'on a 
$a_{\qq}(E)=-6$. D'apr\`es les  conditions (6) et (7), $E_0$ a donc bonne r\'eduction en $\qq$, d'o\`u $v_{\qq}(abc)=0$. On en d\'eduit 
que $89$ divise  $W_{8}$ (cf. d\'em. du lemme 4). Or on a 
$W_{8}=-3^7\times 5^3\times 17^3,$
d'o\`u la contradiction cherch\'ee. On conclut de m\^eme avec les couples $(p,n)$ suivants :
$$(23,20), (53,20), (59,20), (67,4), (79,100), (83,56), (127,4).$$

2)  Compte tenu des \'egalit\'es $W_2=-3$ et $W_{10}=-3\times 11^9\times 31^3$, le corollaire 2 est une cons\'equence directe du th\'eor\`eme.
\bigskip
\medskip
 
\centerline {\douze{Bibliographie}}
\bigskip
\vskip0pt\noindent
[1]  C. Batut, D. Bernardi, K. Belabas, H. Cohen et M. Olivier, PARI-GP, version 2.3.3,  
Universit\'e de Bordeaux I, (2008).
\smallskip
\vskip0pt\noindent
[2] H. Cohen,  Advanced Topics in Computational Number Theory, Springer-Verlag  GTM {\bf{193}}, 2000.
\smallskip
\vskip0pt\noindent
[3]  N. Freitas et S. Siksek, The asymptotic Fermat's Last Theorem for five-sixths of real quadratic fields,
arxiv : 1307.3162v3  (16 Jul 2014), 21 pages, \`a para\^\i tre dans la revue Compositio Mathematica.
\smallskip
\vskip0pt\noindent
[4]  N. Freitas et S. Siksek,   Fermat's Last Theorem over some small real quadratic fields, 
arxiv : 1407.4435v1  (16 Jul 2014), 15 pages.
\smallskip
\vskip0pt\noindent
[5] N. Freitas,  S. Siksek et Le Hung, Elliptic curves over real quadratic fields are mo\-dular, arxiv : 1310.7088v4  (18 Jul 2014), 38 pages, \`a para\^\i tre
dans la revue Inventiones Mathematicae.
\smallskip
\vskip0pt\noindent
[6] H. Gross et D. E. Rohrlich, Some results on the Mordell-Weil group of the Jacobian of the Fermat curve, {\it{Invent. Math.}} {\bf{44}} (1978), 
201-224.
\smallskip
\vskip0pt\noindent
[7] F. H. Hao et  C.  J. Parry, The Fermat equation over quadratic fields, {\it J. Number Theory} {\bf 19} (1984),  115-130.
\smallskip
\vskip0pt\noindent
[8] F. Jarvis et P. Meekin, The Fermat equation over $\Q\bigl(\sqrt{2}\bigr)$, {\it J. Number Theory} {\bf 109} (2004),  182-196.
\smallskip
\vskip0pt\noindent
[9] S. Kamienny et F. Najman, Torsion groups of elliptic curves over quadratic fields,   {\it Acta Arith.} {\bf{152}} (2012), 291-305.
\smallskip
\vskip0pt\noindent
[10] A. Kraus, Sur le d\'efaut de semi-stabilit\'e des courbes elliptiques \`a r\'eduction additive, {\it{Manuscripta Math.}} {\bf{69}} (1990), 353-385.
\smallskip
\vskip0pt\noindent
[11] A. Kraus, Courbes elliptiques semi-stables sur les corps de nombres, {\it{Int. J. Number Theory}} {\bf{3}} (2007), 611-633.
\smallskip
\vskip0pt\noindent
[12]  LMFDB, The $L$-functions and modular forms database, (2014).
\smallskip
\vskip0pt\noindent
[13] I. Papadopoulos, Sur la classification de N\'eron des courbes elliptiques en caract\'e\-ristique r\'esiduelle $2$ et $3$, {\it J. Number Theory} {\bf 44} (1993), 119-152.
\smallskip
\vskip0pt\noindent
[14] P. Ribenboim, Fermat's Last Theorem for Amateurs, Springer-Verlag, 1999.
\smallskip
\vskip0pt\noindent
[15] J.-P. Serre, Propri\'et\'es galoisiennes des points d'ordre fini des courbes elliptiques, {\it{Invent. Math.}} {\bf{15}} (1972), 259-331.
\smallskip
\vskip0pt\noindent
[16] J. H. Silverman, The Arithmetic of Elliptic Curves, Second Edition, Springer-Verlag GTM {\bf{106}}, 2009.
\smallskip
\vskip0pt\noindent
[17] W. Stein et al., A database of elliptic curves over $\Q\bigl(\sqrt{5}\bigr)$-First report (16 pages), arxiv~: 1202.6612v2  (9 Jul 2012), 17 pages, voir Some Relevant Tables, table.pdf.
\smallskip
\vskip0pt\noindent
[18]  A. Wiles, Modular elliptic curves and
Fermat's Last Theorem, {\it {Ann. of Math.}} {\bf{141}} (1995), 443-551.
\bigskip

\line {\hfill{22 septembre 2014}}

\item{} Alain Kraus
\item{}Universit\'e de Paris VI, 
\item{} Institut de Math\'ematiques, 
\item{} 4 Place Jussieu, 75005 Paris,  
\item{} France
\medskip
\vskip0pt\noindent
\item{}e-mail : alain.kraus@imj-prg.fr

 \bye